\documentclass[12pt]{article}
\usepackage{amssymb}

\usepackage{amsmath,latexsym,amsfonts}
\usepackage{graphicx}
\usepackage{srcltx}
\usepackage{amscd}

\newtheorem{lemma}{Lemma}[section]
\newtheorem{coro}[lemma]{Corollary}
\newtheorem{prop}[lemma]{Proposition}
\newtheorem{thm}[lemma]{Theorem}

\makeatletter\@addtoreset{equation}{section}
\renewcommand\theequation{\thesection.\@arabic\c@equation}

\begin{document}
\begin{center}
{\LARGE  Transversal infinitesimal automorphisms on  K\"ahler foliations}

 \renewcommand{\thefootnote}{}
\renewcommand{\thefootnote}{\arabic{footnote}}
\setcounter{footnote}{0} \vspace{0.5cm}

{\large Seoung Dal Jung}
\end{center}
\noindent{\bf Abstract.} 
Let $\mathcal F$ be a K\"ahler foliation on a
compact Riemannian manifold $M$. We study the properties of infinitesimal automorphisms on $(M,\mathcal F)$, and  in particular we concentrate on the transversal conformal field, transversal projective field and transversally holomorphic field. 

\section{Introduction}
Let $(M,\mathcal F)$ be a Riemannian manifold with a Riemannian foliation $\mathcal F$ of codimension $q$. A transversal infinitesimal automorphism on $M$ is an infinitesimal automorphism which preserves the leaves. A transversal infinitesimal automorphism is said to be a transversal Killing field, a transversal conformal field or a transversal projective field if it generates a one parameter family of a transversal infinitesimal isometric, a transversal infinitesimal confomal or a  transversal infinitesimal projective transformation, respectively. Such geometric objects give some important information about the leaf space $M/\mathcal F$. There are several results about infinitesimal automorphisms on Riemannian foliations [\ref{JJ},\ref{JK},\ref{JLK},\ref{KT1},\ref{NT},\ref{Pak}]. Recently, M. J. Jung and S. D. Jung [\ref{JJ}] studied the properties of transversal infinitesimal automorphisms on a compact foliated Riemanian manifold $(M,\mathcal F)$.

In this paper, we investigate the properties of transversal infinitesimal automorphisms on K\"ahler foliations. 
The paper is organized as follows. In Section 2, we review the basic facts on Riemannian foliations.  In Section 3, we review the well-known facts about  infinitesimal automorphisms on Riemannian foliations. In Section 4, we prove that, on K\"ahler foliations, every transversal conformal (or projective) field is a transversal affine field (Theorem 4.2, Theorem 4.5). In particular, if the transversal scalar curvature is a non-zero constant, every transversal conformal field is a transversal Killing field (Theorem 4.3).
In addition, every transversal projective field satisfying some condition is a transversal Killing field (Theorem 4.6). Note that on ordinary manifolds, any affine field is a Killing field, but  on Riemannian foliations, a transversal affine field is not necessarily a transversal Killing field [\ref{JJ}].  
In Section 5, we study  transversally holomorphic fields  and give a vanishing theorem without making the assumption that all leaves of $\mathcal F$ are minimal.
\section{Preliminaries}
 Let $(M,g_M,\mathcal F)$ be a $(p+q)$-dimensional
Riemannian manifold with a foliation $\mathcal F$ of codimension
$q$ and a bundle-like metric $g_M$ with respect to $\mathcal F$
[16]. Let $\nabla^M$ be the Levi-Civita connection with respect to
$g_M$. Let $TM$ be the tangent bundle of $M$, $L$ its integrable
subbundle given by $\mathcal F$, and  $Q=TM/L$ the corresponding
normal bundle. Then there exists an exact sequence of vector
bundles
\begin{equation}
 0 \longrightarrow L \longrightarrow
TM_{\buildrel \longleftarrow \over \sigma }^{\buildrel \pi \over \longrightarrow} Q \longrightarrow
0,
\end{equation}
where $\sigma:Q\to L^\perp$ is a bundle map satisfying
$\pi\circ\sigma=id$.  Let $g_Q$ be the holonomy invariant metric
on $Q$ induced by $g_M=g_L + g_{L^\perp}$; that is,
\begin{equation}
g_Q(s,t)=g_M(\sigma(s),\sigma(t)) \quad \forall\ s,t\in \Gamma Q.
\end{equation}
This means that $\theta(X)g_Q=0$ for $X\in\Gamma L$, where
$\theta(X)$ is the transverse Lie derivative. So we have an
identification $L^\perp$ with $Q$ via an isometric splitting
$(Q,g_Q)\cong (L^\perp, g_{L^\perp})$. A transversal Levi-Civita
connection $\nabla$ in $Q$ is defined [5] by
\begin{equation}\label{connection}
\nabla_X s=\left\{\begin{split}& \pi([X,Y_s])\qquad \forall
X\in\Gamma L\\
 &\pi(\nabla^M_X Y_s)\qquad \forall X\in\Gamma L^\perp,
\end{split}
\right.
\end{equation}
where $s\in \Gamma Q$ and $Y_s =\sigma(s)\in \Gamma L^\perp $
corresponding to $s$ under the canonical isomorphism $Q \cong
L^\perp$. The curvature $R^\nabla$ of $\nabla$ is defined by $
R^\nabla(X,Y)=[\nabla_X,\nabla_Y]-\nabla_{[X,Y]}$  for
$X,Y\in\Gamma TM$. Since $i(X)R^\nabla=0$ for any $X\in \Gamma L$
[5], we can define the transversal Ricci operator
$\rho^\nabla:\Gamma Q \to \Gamma Q$ by
\begin{equation}
\rho^\nabla (s_x)=\sum_{a=p+1}^n R^\nabla(s,E_a)E_a,
\end{equation}
where $\{E_a\}_{a=p+1,\cdots,n}$ is an orthonormal basic frame of
$Q$. And then the transversal Ricci curvature ${\rm Ric^\nabla}$
is given by
$\operatorname{Ric^\nabla}(s_1,s_2)=g_Q(\rho^\nabla(s_1),s_2)$ for
any $s_1,s_2\in\Gamma Q$. The transversal scalar curvature
$\sigma^\nabla$ is given by $\sigma^\nabla = {\rm Tr}
\rho^\nabla$.
 The foliation $\mathcal F$ is said to be (transversally) {\it Einsteinian}  if the model space
 is Einsteinian, that is,
\begin{equation}
\rho^\nabla=\frac1q \sigma^\nabla\cdot {\rm id}
\end{equation}
with constant transversal scalar curvature  $\sigma^\nabla.$
 The {\it mean curvature vector} $\kappa^\sharp$ of $\mathcal F$ is defined by
\begin{equation}
\kappa^\sharp = \pi \big(\sum_{i=1}^p \nabla_{E_i}^M E_i \big),
\end{equation}
where $\{E_i\}$ is a local orthonormal basis of $L$. The foliation $\mathcal F$ is said to be {\it
minimal} if $\kappa^\sharp=0$.
A differential form $\omega\in \Omega^r(M)$ is {\it basic} if $
i(X)\omega=0$ and $\theta(X)\omega=0$ for all $X\in\Gamma L$.
Let $\Omega_B^r(\mathcal F)$ be the set of all basic r-forms on
$M$. Then $\Omega^r(M)=\Omega_B^r(\mathcal F)\oplus \Omega_B^r(\mathcal F)^\perp$ [\ref{L}]. It is well-known that the mean curvature form
$\kappa_B$ is closed, i.e., $d\kappa_B=0$, where $\kappa_B$ is the basic part of $\kappa$.
   The {\it basic Laplacian} acting on
$\Omega_B^*(\mathcal F)$ is defined by
\begin{equation}
\Delta_B=d_B\delta_B+\delta_B d_B,
\end{equation}
where $\delta_B$ is the formal adjoint of
$d_B=d|_{\Omega_B^*(\mathcal F)}$ [\ref{L},\ref{Jung}]. 
Let
 $\{E_a\}(a=1,\cdots,q)$ be a local orthonormal basis of $Q$. We define  $\nabla_{\rm tr}^*\nabla_{\rm tr}:\Omega_B^r(\mathcal F)\to \Omega_B^r(\mathcal F)$ by
\begin{align}\label{eq1-12}
\nabla_{\rm tr}^*\nabla_{\rm tr} =-\sum_a \nabla^2_{E_a,E_a}
+\nabla_{\kappa_B^\sharp},
\end{align}
where $\nabla^2_{X,Y}=\nabla_X\nabla_Y -\nabla_{\nabla^M_XY}$ for
any $X,Y\in TM$. The operator $\nabla_{\rm tr}^*\nabla_{\rm tr}$
is positive definite and formally self adjoint on the space of
basic forms [\ref{Jung}]. We define the bundle map $A_Y:\Lambda^r
Q^*\to\Lambda^r Q^*$ for any $Y\in V(\mathcal F)$ [\ref{KT1}]
by
\begin{align}\label{eq1-13}
A_Y\phi =\theta(Y)\phi-\nabla_Y\phi,
\end{align}
where $\theta(Y)$ is the transverse Lie derivative. 
Then it is proved [6] that, for any vector field
$Y\in V(\mathcal F)$,
\begin{equation}\label{2-4}
A_Y s = -\nabla_{Y_s}\bar Y,
\end{equation}
where $Y_s =\sigma(s)\in \Gamma TM$.   So $A_Y$ depends only on
$\bar Y=\pi(Y)$ and is a linear operator.
Since
$\theta(X)\phi=\nabla_X\phi$ for any $X\in\Gamma L$, $A_Y$
preserves the basic forms and depends only on $\bar Y$.    
     Then we
have the generalized Weitzenb\"ock formula.
\begin{thm} $[\ref{Jung}]$ On a Riemannian foliation $\mathcal F$, we have
\begin{align}\label{eq2-11}
  \Delta_B \phi = \nabla_{\rm tr}^*\nabla_{\rm tr}\phi +
  F(\phi)+A_{\kappa_B^\sharp}\phi,\quad\phi\in\Omega_B^r(\mathcal
  F),
\end{align}
 where $F(\phi)=\sum_{a,b}\theta^a \wedge i(E_b)R^\nabla(E_b,
 E_a)\phi$. If $\phi$ is a basic 1-form, then $F(\phi)^\sharp
 =\rho^\nabla(\phi^\sharp)$.
\end{thm}
From Theorem 2.1, we have the following. For any $\phi\in\Omega_B^r(\mathcal F)$,
\begin{align}\label{eq2-12}
\frac12\Delta_B|\phi|^2 = \langle\Delta_B\phi,\phi\rangle -|\nabla_{\rm tr}\phi|^2 -\langle F(\phi),\phi\rangle -\langle A_{\kappa_B^\sharp}\phi,\phi\rangle.
\end{align}
Now, we recall the following generalized maximum principle.
\begin{lemma} $[\ref{JLK}]$ Let  $\mathcal F$ be a Riemannian foliation on a closed, oriented Riemannian manifold $(M,g_M)$. If $(\Delta_B -\kappa_B^\sharp)f\geq 0$ $($or $\leq 0)$ for any basic function $f$, then $f$ is constant.
\end{lemma}
 Let $V(\mathcal F)$ be the space of all vector fields $Y$ on $M$ satisfying $[Y,Z]\in \Gamma L$ for
all $Z\in \Gamma L$. An element of $V(\mathcal F)$ is called an
{\it infinitesimal automorphism} of $\mathcal F$ [\ref{NT}]. Let
\begin{equation}\label{2-7}
\bar V(\mathcal F)=\{\bar Y=\pi(Y)| Y\in V(\mathcal F)\}.
\end{equation}
It is trivial that an element $s$ of $\bar V(\mathcal F)$
satisfies $\nabla_X s=0$ for all $X\in \Gamma L$ [\ref{KT1}]. Hence $
\bar V(\mathcal F)\cong \Omega_B^1(\mathcal F)$.

\section{Transversal infinitesimal automorphisms}
 If $Y\in V(\mathcal F)$ satisfies $\theta(Y)g_Q=0$, then $\bar Y$ is called a {\it transversal
Killing field} of $\mathcal F$.
 If $Y\in V(\mathcal F)$ satisfies $\theta(Y)g_Q=2f_Y g_Q$ for
a basic function $f_Y$ depending on $Y$, then $\bar Y$ is called a
{\it transversal conformal field} of $\mathcal F$. Equivalently, for any $X,Z\in V(\mathcal F)$
\begin{align}
g_Q(\nabla_X\bar Y,Z)+g_Q(X,\nabla_Z\bar Y)=2f_Y g_Q(\bar X,\bar Z).
\end{align}
In this case,
we have
\begin{equation}
f_Y = \frac 1q \operatorname{div_\nabla} \bar Y,
\end{equation}
where ${\rm div}_\nabla \bar Y$ is the transversal divergence of $\bar Y$.
A transversal conformal field $\bar Y$ is {\it homothetic} if $f_Y$ is constant. 
 For any vector fields $Y,Z\in V(\mathcal F)$ and $X\in\Gamma Q$, we have [\ref{JJ}]
\begin{equation}\label{eq3-2}
(\theta(Y)\nabla)(Z,X)=R^\nabla(\bar Y,\bar Z)X +\nabla_{\bar Z}\nabla_X \bar Y
-\nabla_{\nabla_ZX}\bar Y.
\end{equation}
 If $Y\in V(\mathcal F)$ satisfies $\theta(Y)\nabla=0$, then $\bar Y$ is called a {\it transversal
 affine field} of $\mathcal F$.  If $Y\in V(\mathcal F)$ satisfies
 \begin{equation}\label{eq3-3}
( \theta(Y)\nabla)(X,Z)=\alpha_Y(X)Z +\alpha_Y(Z)X
\end{equation}
for any $X,Z\in\Gamma Q$, where $\alpha_Y$ is a basic 1-form on $M$, then $\bar Y$ is called
a {\it transversal projective field} of $\mathcal F$; in this
case, it is trivial that
\begin{align}\label{eq3-4}
\alpha_Y = {1\over q+1} d_B \operatorname{div_\nabla} \bar Y.
\end{align} 
Let $\{E_a\}_{a=1,\cdots,q}$ be a local orthonormal basic  frame in $Q$ such that $(\nabla E_a)_x$ for $x\in M$. From now on, all the computations in this paper will be made in such charts.
For any $Y\in V(\mathcal F)$, from (\ref{eq3-2}), we have
\begin{align}\label{eq3-5}
(\theta(Y)R^\nabla)(E_a,E_b)E_c =(\nabla_a\theta(Y)\nabla)(E_b,E_c)-(\nabla_b\theta(Y)\nabla)(E_a,E_c),
\end{align}
where $\nabla_a = \nabla_{E_a}$. Then we have the following lemma.
 \begin{lemma} $[\ref{JJ}]$  Let $\mathcal F$ be a Riemannian foliation of codimension $q$ on a Riemannian manifold $(M,g_M)$. If $\bar Y \in
\bar V(\mathcal F)$ is a transversal conformal  field, i.e.,
$\theta(Y)g_Q=2f_Y g_Q$, then we have
\begin{align}
&g_Q((\theta(Y) \nabla)(E_a,E_b),E_c)=\delta_b^c f_a + \delta_a^c
f_b -\delta_a^b
f_c,\label{eq3-5}\\
&g_Q((\theta(Y) R^\nabla)(E_a,E_b)E_c,E_d)=\delta_b^d \nabla_a f_c
- \delta_b^c \nabla_a f_d
-\delta_a^d \nabla_b f_c +\delta_a^c \nabla_b f_d, \label{eq3-6}\\
&\theta(Y)\sigma^\nabla=2(q-1)(\Delta_B f_Y
-\kappa_B^\sharp(f_Y))-2f_Y\sigma^\nabla,\label{eq3-7}
\end{align}
where  $f_a=\nabla_a f_Y$.
\end{lemma}
From (\ref{eq3-5}), it is trivial that any transversal homothetic field is a transversal affine field.
On the other hand, from (\ref{eq3-3}) and (\ref{eq3-5}), we have the following.
\begin{lemma} Let $\mathcal F$ be the same as in Lemma 3.1. If $\bar Y \in
 \bar V(\mathcal F)$ is a transversal projective  field, then we have
\begin{align}
(\theta(Y)R^\nabla)(E_a,E_b)E_c&=(\nabla_a\alpha_Y)(E_b)E_c +(\nabla_a\alpha_Y)(E_c)E_b\\
         &-(\nabla_b \alpha_Y)(E_a)E_c -(\nabla_b\alpha_Y)(E_c)E_a.\notag
\end{align}
\end{lemma}
Now, we define the operator $B^\mu_Y :\Gamma Q\to \Gamma Q (\mu\in\mathbb R)$ for any $Y\in V(\mathcal F)$ by
\begin{align}
B_Y=A_Y + A_Y^t+\mu\cdot {\rm div}_\nabla \bar Y {\rm id}.
\end{align}
It is well-known [\ref{KT1}] that $\bar Y$ is transversal conformal(resp. transversal Killing) field if and only if $B_Y^{2/q}=0\ ({\rm resp.}\ B_Y^0=0)$.
Then we have the following lemma.
\begin{lemma} Let $\mathcal F$ be a Riemannian foliation  on a compact Riemannian manifold $(M,g_M)$. If $\bar Y$ is transversal homothetic, i.e., ${\rm div}_\nabla\bar Y$ is constant, then
\begin{align}
\int_M g_Q(B_Y^\mu \bar Y,\kappa^\sharp_B)={\rm div}_\nabla\bar Y(\mu-{2\over q}{\rm div}_\nabla\bar Y) {\rm vol}(M).
\end{align}
\end{lemma}
{\bf Proof.} From (3.1) and the transversal divergence theorem, the proof follows. $\Box$

Now, we recall the following relationships among infinitesimal automorphisms on a Riemannian foliation.
\begin{prop} $[\ref{JJ}]$ Let $\mathcal F$ be a Riemannian foliation  on a compact Riemannian manifold $(M,g_M)$. Then,

 $(1)$ Any transversal Killing field is a transversal affine field.

 $(2)$ Any transversal affine field with $\int_M g_Q(B_Y^0\bar Y,\kappa^\sharp)=0$ is a transversal Killing field.

 $(3)$ Any transversal conformal field (or projective field) $\bar Y$ with the properties
 \begin{align*}
   (i) \ \int_M g_Q(B_Y^0 \bar Y,\kappa^\sharp)\geq 0,\quad  (ii)\ d_B {\rm div}_\nabla \bar Y=0
   \end{align*}
    is a transversal Killing field.
\end{prop}
Note that on $\mathcal F$ with a constant transversal scalar curvature $\sigma^\nabla$, if $\mathcal F$ admits a transversal conformal field $\bar Y$ with $f_Y\ne 0$, then $\sigma^\nabla$ is non-negative ([\ref{JJ}], Corollary 5.6). Equivalently, on $\mathcal F$ with a negative constant $\sigma^\nabla$, there is no non-isometric transversal conformal field.
Hence we have the following proposition.
\begin{prop} Let $\mathcal F$ be Riemannian foliation of codimension $q$ on a compact Riemannian manifold $(M,g_M)$. Assume that the transversal scalar curvature $\sigma^\nabla$ is negative constant. Then any transversal conformal field is a transversal Killing field.
\end{prop}
\begin{thm} $[\ref{JK}]$ Let $(M,g_M,\mathcal F)$ be a compact Riemannian manifold with a foliation $\mathcal F$ of codimension $q$ and a bundle-like metric $g_M$.  Assume that the transversal Ricci curvature $\rho^\nabla$ is non-positive and negative at some point. Then

$(1)$ there are no transversal Killing fields on $M$ 

$(2)$ if $\delta_B\kappa_B=0$, then there are no transversal conformal fields.
\end{thm} 
\noindent{\bf Remark.} From Proposition 3.5 and Theorem 3.6, it is well-known that on a transversally Einstein foliation with negative scalar curvature, there are no transversal conformal fields without the condition $\delta_B\kappa_B=0$. For more relations among infinitesimal automorphisms on a Riemannian foliation, see [\ref{JJ},\ref{Pak}].
\section{Transversal conformal and projective field on K\"ahler foliations}
Now, we study the infinitesimal automorphisms on K\"ahler foliations. 
Let $\mathcal F$ be a K\"ahler foliation of codimension $q=2m$ on a Riemannian manifold $(M,g_M)$ [\ref{NT}].
 Note that for any $X,
Y\in\Gamma Q$,
\begin{equation}
\Omega(X,Y)=g_Q(X,JY)
\end{equation}
defines a basic 2-form $\Omega$, which is closed, where $J:Q\to Q$ is an almost
complex structure on $Q$.  Then we have
\begin{align}
\Omega=-\frac12\sum_{a=1}^{2m}\theta^a\wedge J\theta^a,
\end{align}
where $\theta^a$ is a dual form of $E_a$.
  Moreover, we have the following
identities:
\begin{gather}\label{eq3-10}
 R^{\nabla}(X,Y)J = JR^{\nabla}(X,Y),\quad
 R^{\nabla}(JX,JY) =
R^{\nabla}(X,Y)
\end{gather}
for any  $X, Y \in\Gamma Q$. Then we have the following.
\begin{prop} \label{prop3-5}  Let $\mathcal F$ be a K\"ahler foliation of codimension $q=2m$ on a  Riemannian manifold $(M,g_M)$
 and let $\bar Y$ be a transversal conformal field, i.e.,
 $\theta(Y)g_Q=2f_Y g_Q$. Then we have
 \begin{align}\label{eq3-17}
 \Delta_B f_Y -\kappa_B^\sharp(f_Y)=0.
 \end{align}
 Moreover, if $M$ is compact, then $f_Y$ is constant, i.e., $\bar Y$ is transversal homothetic.
 \end{prop}
 {\bf Proof.} Let $f_Y$ be a basic function with $\theta(Y)g_Q=f_Y
 g_Q$. Fix $x\in M$ and let $\{E_a\}$ be a local orthonormal basic frame such that $(\nabla E_a)_x=0$. 
Then, at $x$, we have from (\ref{eq3-6}), 
\begin{align*}
&\sum_{a,b=1}^{2m}g_Q((\theta(Y)R^\nabla)(E_a,E_b)E_a,E_b)=2q\sum_{a=1}^{2m} E_a
E_a (f_Y)\\
&\sum_{a,b=1}^{2m}g_Q((\theta(Y)R^\nabla)(JE_a,JE_b)E_a,E_b)=2\sum_{a=1}^{2m} E_a
E_a (f_Y).
\end{align*}
From (\ref{eq3-10}),  we have
\begin{align}\label{eq3-18}
2(q-1)\sum_{a=1}^{2m} E_a E_a (f_Y)=0.
\end{align}
Since $q>1$, $\sum_{a=1}^{2m} E_aE_a(f_Y)=0$. Hence  $\Delta_B f_Y =\kappa_B^\sharp(f_Y)$, which proves (\ref{eq3-17}). Moreover, if $M$ is compact, by Lemma 2.2, $f_Y$ is constant. $\Box$
\begin{thm} Let $\mathcal F$ be a K\"ahler foliation on a compact Riemannian manifold
$(M,g_M)$. Then any transversal conformal field is a tranversal affine field.
\end{thm}
{\bf Proof.} Let  $\bar Y$ be a transversal conformal field
such that $\theta(Y)g_Q=f_Y g_Q$. By Proposition \ref{prop3-5}, $f_Y$ is constant. Therefore,  from Lemma 3.1 (\ref{eq3-5}),  $\theta(Y)\nabla=0$. So $\bar Y$ is the transversal affine field. $\Box$

\noindent{\bf Remark.} On a compact K\"ahler manifold, any conformal field is always a Killing field, because any affine field is a Killing field [\ref{TA}]. For the foliated manifold, this does not hold because of Proposition 3.3 (2). In fact, we have the following theorem.

\begin{thm} Let $\mathcal F$ be a K\"ahler foliation on a compact Riemannian manifold $(M,g_M)$. Assume that the transversal scalar
curvature $\sigma^\nabla$ is non-zero constant. Then any
transversal conformal field is a transversal Killing
field.
\end{thm}
{\bf Proof.} Let $\bar Y$ be a transversal conformal  field
such that $\theta(Y)g_Q=f_Y g_Q$.  Since $\sigma^\nabla\ne 0$ is
constant, from Lemma 3.1 (\ref{eq3-7}) and Proposition \ref{prop3-5}, $f_Y=0$.
Therefore, $\bar Y$ is a transversal Killing  field.
$\Box$

Now we study the transversal projective field on a K\"ahler foliation. From Lemma 3.2, we have the following.
\begin{prop} Let $\mathcal F$ be a K\"ahler foliation of codimension $q=2m (m\geq2)$ on a Riemannian manifold $(M,g_M)$ and let $\bar Y$  be a transversal projective field. Then we have
\begin{align}\label{eq3-14}
\Delta_B g_Y -\kappa_B^\sharp(g_Y)=0,
\end{align}
where $g_Y= {\rm div}_\nabla(\bar Y)$. If $M$ is compact, then $g_Y$ is constant.
\end{prop}
{\bf Proof.} Let $\bar Y$ be a transversal projective field. Let $\{E_a\}$ be a local orthonormal basic frame such that $(\nabla E_a)_x=0$ at $x\in M$. Then, from Lemma 3.2, we have
\begin{align*}
(q+1)(\theta(Y)R^\nabla)(E_a,E_b)E_c &= E_a E_b({\rm div}_\nabla\bar Y) E_c + E_a E_c ({\rm div}_\nabla\bar Y)E_b\\
&-E_b E_a({\rm div}_\nabla\bar Y)E_c -E_b E_c({\rm div}_\nabla\bar Y)E_a.
\end{align*}
Hence we have
\begin{align*}
&(q+1)\sum_{a,b=1}^{2m}g_Q((\theta(Y)R^\nabla)(E_a,E_b)E_a,E_b)=(q-1)\sum_{a=1}^{2m} E_a
E_a (f_Y)\\
&(q+1)\sum_{a,b=1}^{2m}g_Q((\theta(Y)R^\nabla)(JE_a,JE_b)E_a,E_b)=\sum_{a=1}^{2m} E_a
E_a (f_Y).
\end{align*}
From (\ref{eq3-10}), we have
\begin{align*}
(q-2)\sum_{a=1}^{2m} E_a E_a ({\rm div}_\nabla \bar Y) =0.
\end{align*}
Since $q>2$, we have $(\Delta_B-\kappa_B^\sharp)g_Y =\sum_{a=1}^{2m} E_a E_a (g_Y)=0$. Moreover, if $M$ is compact, by Lemma 2.2, $g_Y$ is constant. $\Box$

From Proposition 4.4, we have the following theorem.
\begin{thm}  Let $\mathcal F$ be a K\"ahler foliation of codimension $q=2m(m\geq2)$ on a compact Riemannian manifold
$(M,g_M)$. Then any transversal projective field is a tranversal
affine field.
\end{thm} 
{\bf Proof.} Let $\bar Y$ be a transversal projective field. From Proposition 4.4,  $g_Y={\rm div}_\nabla\bar Y$ is constant. Hence, from (\ref{eq3-4}), $\alpha_Y=0$. From (\ref{eq3-3}), $\bar Y$ is transversal affine. $\Box$

From Proposition 3.3 (3), we have the following corollary.
\begin{coro}  Let $\mathcal F$ be a K\"ahler foliation of codimension $q=2m(m\geq2)$ on a compact Riemannian manifold
$(M,g_M)$. If any transversal projective field $\bar Y$ satisfies $\int_M g_Q(B_Y^0 \bar Y,\kappa^\sharp)\geq 0$, then $\bar Y$ is a transversal Killing field.
\end{coro}
{\bf Remark.} For a harmonc K\"ahler foliation on a compact Riemannian manifold, any transversal projective is a transversal Killing field. For the point foliation, any transversal affine field is a transversal Killing field [\ref{Yano}]. So Theorem 4.2 and Theorem 4.5 are found in [\ref{TA}].

\section{Transversally holomorphic fields}
Let $\mathcal F$ be a K\"ahler foliation of codimension $q=2m$ on a Riemannian manifold $(M,g_M)$. Let $Y$ be an infinitesimal automorphism of $\mathcal F$. Then a vector field $\bar Y$ is said to be a {\it transversally holomorphic field}  if
\begin{align}
\theta(Y)J=0,
\end{align}
equivalently,  for all $Z\in\Gamma L^\perp$
\begin{align}
\nabla_{JZ}\bar Y = J\nabla_Z \bar Y.
\end{align}
Let $\{E_\alpha,JE_\alpha\} (\alpha=1,\cdots,m)$ be a local orthonormal basis of $\Gamma L^\perp$.
Then we recall the following well-known facts.
\begin{lemma} $[\ref{NT}]$ On a K\"ahler foliation of codimension $q=2m$, it holds that
\begin{align}\label{5-3}
\rho^\nabla(X)=\sum_{a=1}^m J R^\nabla(E_a,JE_a)X.
\end{align}
\end{lemma}
Then we have the following facts on a harmonic foliation.
\begin{thm} $[\ref{NT}]$ On a harmonic K\"ahler foliation $\mathcal F$ on a compact manifold $(M,g_M)$, 
the followings are equivalent:

$(1)$ $\bar Y$ is transversally holomorphic, $\theta (Y) J=0$,

$(2)$ $\bar Y$ is a transversal Jacobi field of $\mathcal F$, i.e., $\nabla_{\rm tr}^*\nabla_{\rm tr}\bar Y-\rho^\nabla(\bar Y)=0$.
\end{thm}
Next, we study the above relations on a non-harmonic K\"ahler foliation. In fact, we have the following theorem.
\begin{thm} On a K\"ahler foliation $\mathcal F$ on $(M,g_M)$,  $\bar Y$ is transversally holomorphic, i.e.,  $\theta(Y)J=0$ if and only if
\begin{align*}
&(i)\ \nabla_{\rm tr}^*\nabla_{\rm tr}\bar Y-\rho^\nabla(\bar Y)+A_Y\kappa_B^\sharp=0, \\
&(ii)\ \int_Mg_Q((\theta(Y)J)\kappa_B^\sharp,J\bar Y)=0.
\end{align*}
\end{thm}
{\bf Proof.} Let $\bar Y$ be transversally holomorphic, i.e., $\nabla_JZ \bar Y=J\nabla_Z\bar Y$ for any $Z\in \Gamma Q$. Then, by long calculation, we have
\begin{align}\label{5-4}
\nabla_{\rm tr}^*\nabla_{\rm tr}\bar Y&=\sum_{\alpha=1}^m JR^\nabla(E_\alpha,JE_\alpha)\bar Y +\nabla_{\kappa_B^\sharp}\bar Y.
\end{align}
From (\ref{2-4}) and (\ref{5-3}), we have
\begin{align}
 \nabla_{\rm tr}^*\nabla_{\rm tr}\bar Y-\rho^\nabla(\bar Y)+A_Y\kappa_B^\sharp=0.
 \end{align}
 Hence (i) and (ii) are proved.  Conversely, by direct calculation, we have
 \begin{align*}
 \int_M |\theta(Y)J|^2 &=2\int_M g_Q(\nabla_{\rm tr}^*\nabla_{\rm tr}\bar Y-\rho^\nabla(\bar Y)+A_Y\kappa_B^\sharp,\bar Y) \\
 &+2\int_M\sum_{a=1}^{2m} E_ag_Q(\nabla_{E_a}\bar Y + J\nabla_{JE_a}\bar Y,\bar Y).
 \end{align*}
 Now we choose $X\in \Gamma Q$ by $g_Q( X,Z) =g_Q( \nabla_Z\bar Y+J\nabla_{JZ}\bar Y,\bar Y)$ for any $Z\in\Gamma Q$. Then, by the transversal divergence theorem, we have
 \begin{align*}
 \int_M\sum_{a=1}^{2m} E_ag_Q( \nabla_{E_a}\bar Y+J\nabla_{JE_a}\bar Y,\bar Y) &=\int_M div_\nabla (X)\\
&=\int_M g_Q(\nabla_{\kappa_B^\sharp}\bar Y+ J\nabla_{J\kappa_B^\sharp}\bar Y,\bar Y).
\end{align*}
Hence we have 
 \begin{align*}
\frac12 \int_M |\theta(Y)J|^2 &=\int_M g_Q( \nabla_{\rm tr}^*\nabla_{\rm tr}\bar Y-\rho^\nabla(\bar Y)+A_Y\kappa_B^\sharp,\bar Y) +\int_M g_Q((\theta(Y)J)\kappa_B^\sharp,J\bar Y).
 \end{align*}
 Hence the converse is proved. $\Box$.
 
 \noindent{\bf Remark.} The solution of $ \nabla_{\rm tr}^*\nabla_{\rm tr}\bar Y-\rho^\nabla(\bar Y)+A_Y\kappa_B^\sharp=0$ appears as the kernel of the transversal Jacobi operator $J_{\rm id}^T$ of the identity map [\ref{Jung1}]. 

Moreover, on non-harmonic K\"ahler foliations on compact manifolds, the following theorem holds.
\begin{thm} Let $\mathcal F$ be a K\"ahler foliation on a compact Riemannian manifold $(M,g_M)$. Assume that the transversal Ricci operator is non-positive and negative at some point. Then every infinitesimal automorphism $Y$ with a transversally holomorphic field $\bar Y$ satisfies $Y\in\Gamma L$, i.e., $\bar Y=0$.
 \end{thm}
 {\bf Proof.} Let $\bar Y$ be a transversally holomorphic field. Then, by Theorem 5.3, we have
 \begin{align*}
 \Delta_B |\bar Y|^2 &= 2 g_Q(\nabla_{\rm tr}^*\nabla_{\rm tr}\bar Y,\bar Y) -2|\nabla_{\rm tr}\bar Y|^2\\
  &=2g_Q(\rho^\nabla(\bar Y),\bar Y)-2|\nabla_{\rm tr}\bar Y|^2 +\kappa_B^\sharp|\bar Y|^2.
  \end{align*}
  Since the transversal Ricci curvature $\rho^\nabla$ is non-positive, we have $(\Delta_B-\kappa_B^\sharp)|\bar Y|^2\leq 0$. Hence, by Lemma 2.2, $|\bar Y|$ is constant. Moreover, since $\rho^\nabla$ is negative at some point, $\bar Y$ is zero, i.e., $Y$ is tangential to $\mathcal F$. $\Box$
   
\noindent{\bf Remark.} In [\ref{NT}], S. Nishikawa and P. Tondeur proved Theorem 5.4 when the foliation is minimal.

\bigskip
\noindent{\bf Acknowledgements} This research was supported by the Basic Science Research Program through the National Research Foundation of Korea(NRF) funded by the Ministry of Education, Science and Technology (2010-0021005)

\noindent Department of Mathematics and Research Institute for Basic Sciences, Jeju National University,
Jeju 690-756, Korea

\noindent e-mail: sdjung@jejunu.ac.kr

\end{document}